\documentclass[11pt]{amsart}
\usepackage{amssymb}

\theoremstyle{plain}
\newtheorem{theorem}{Theorem}

\newtheorem{lemma}{Lemma}
\newtheorem{proposition}{Proposition}

\theoremstyle{definition}

\theoremstyle{remark}

\numberwithin{equation}{section}
\begin{document}
\title[Symmetry, Specializability and Infinite Products. ]
       {Symmetry and Specializability in the continued fraction
expansions
of some infinite products}
\author{J. Mc Laughlin}
\address{
}
\email{jgmclaug@math.uiuc.edu}
\begin{abstract}
Let $f(x) \in \mathbb{Z}[x]$. Set $f_{0}(x) = x$ and, for $n \geq 1$, define  $f_{n}(x)$ $=$
$f(f_{n-1}(x))$.
We describe several infinite families of polynomials for which the infinite product
{\allowdisplaybreaks
\begin{equation*}
 \prod_{n=0}^{\infty}\left ( 1 + \frac{1}{f_{n}(x)} \right )
\end{equation*}
}
has a \emph{specializable} continued fraction expansion of the form
{\allowdisplaybreaks
\begin{equation*}
S_{\infty} = [1;a_{1}(x),\,a_{2}(x),\, a_{3}(x), \dots ],
\end{equation*}
}
where $a_{i}(x) \in \mathbb{Z}[x]$, for $i \geq 1$.

When the infinite product and the continued fraction are
\emph{specialized} by letting $x$ take integral values, we get
infinite classes of real numbers whose regular continued fraction
expansion is predictable.

We also show that, under some simple conditions,
all the real numbers produced by this specialization are transcendental.
\end{abstract}

\maketitle


\section{Introduction } \label{S:intro}
The problem of finding the regular continued fraction expansion of
an irrational quantity expressed in some other form has a long
history but until the 1970's not many examples of such continued
fraction expansions were known. Apart from the quadratic irrationals and
numbers like $e^{q}$, for certain rational $q$, there were very
few examples of irrational numbers with predictable patterns in
their sequence of partial quotients.

Being able to predict a pattern in the regular continued fraction expansion
of an irrational number is not only interesting in its own right,
but if one can also derive sufficient  information about the convergents, it is
then sometimes  possible to prove that the number is
transcendental.

In \cite{L73}, Lehmer showed that certain quotients of modified
Bessel functions evaluated at various rationals had  continued
fraction expansions in which the partial quotients lay in
arithmetic progressions. He also showed that similar quotients of
modified Bessel functions evaluated at the square root of a
positive integer had  continued fraction expansions in which the
sequence of partial quotients consisted of  interlaced
arithmetic progressions.

An old result, originally due to B\"{o}hmer \cite{B27} and
Mahler \cite{M29},
was rediscovered by Davison \cite{D77} and Adams and Davison
 \cite{AD77} (generalizing Davison's previous result
in \cite{D77}).
In this latter paper, the authors were able
to determine, for any positive integer $a \geq 2$ and any
positive irrational number $\alpha$,
 the regular continued fraction expansion of the number
{\allowdisplaybreaks
\begin{align}\label{ex1}
S_{a}(\alpha)=(a-1)\sum_{r=1}^{\infty}\frac{1}{a^{\lfloor \alpha \rfloor}}
\end{align}
}
in terms of the convergents in the continued fraction expansion
of $\alpha^{-1}$. They were further able to show that all such
numbers $S_{a}(\alpha)$ are transcendental.

A generalization of Davison's result from \cite{D77} was given by
Bowman in \cite{B88} and  Borwein and Borwein \cite{BB93} gave a
two-variable generalization of \eqref{ex1}  but  the continued
fraction expansion in this latter  case is not usually regular.

Shallit \cite{S79} and      Kmo\u{s}ek \cite{K79} showed
independently  that the continued fraction expansions
of the irrational numbers
{\allowdisplaybreaks
\begin{align*}
\sum_{k=0}^{\infty}\frac{1}{u^{2^{k}}},
\end{align*}
} have predictable continued fraction expansions. This result was
subsequently generalized by K\"{o}hler \cite{K80}, Peth\"{o}
\cite{P82} and by Shallit \cite{S82} once again.

In \cite{MVDPS91}, Mend\`{e}s and van der Poorten considered infinite
products of the form
{\allowdisplaybreaks
\begin{align*}
\prod_{h=0}^{\infty}\left(1+X^{-\lambda_{h}}\right),
\end{align*}
} where $0 < \lambda_{1}<\lambda_{2}< \cdots $ is any sequence of rational integers
satisfying a certain growth condition and showed that such
products had a predictable continued fraction expansion in which
all the partial quotients were polynomials  in
$\mathbb{Z}[\,X\,]$. They further showed that if the infinite
product and continued fraction were \emph{specialized} by letting
$X$ be any integer $g\geq 2$, that all such real numbers
{\allowdisplaybreaks
\begin{align*}
\gamma = \prod_{h=0}^{\infty}\left(1+g^{-\lambda_{h}}\right),
\end{align*}
}
so obtained were transcendental.
Similar investigations, in which the continued fraction expansions
of certain formal Laurent series are determined, can be found in
\cite{VDPS92}, \cite{SVDPS93}, \cite{VDP93} and \cite{ALMSV96}.

Let $f(x)\in \mathbb{Z}[\,x\,]$, $f_{0}(x)=x$ and, for $i \geq 1$,
$f_{i}(x)=f(f_{i-1}(x))$, the $i$-th iterate of $f(x)$. In
\cite{T91}, Tamura investigated infinite series of the
form
{\allowdisplaybreaks
\begin{align*}
\theta(x:f)=\sum_{m=0}^{\infty}\frac{1}{f_{0}(x)f_{1}(x)\cdots f_{m}(x)}.
\end{align*}
}
He showed that, for all polynomials in a certain congruence class,
 that the continued fraction expansion of $\theta(x:f)$ had
all partial quotients in $\mathbb{Z}[\,x\,]$. He further showed that
if the series and continued fraction were specialized to
a sufficiently large integer (depending on $f(x)$) that the resulting
number was transcendental.

The infinite series $\sum_{k=0}^{\infty}1/x^{2^{k}}$, investigated by
Shallit \cite{S79} and      Kmo\u{s}ek \cite{K79}
may be regarded as a special case of the infinite series
$\sum_{k=0}^{\infty}1/f_{k}(x)$, with $f(x) = x^{2}$. In a
very interesting paper, \cite{C96}, Cohn gave a complete classification
of all those polynomials $f(x) \in \mathbb{Z}[\,x\,]$ for which
the series $\sum_{k=0}^{\infty}1/f_{k}(x)$ had a continued fraction
expansion in which all partial quotients were in $\mathbb{Z}[\,x\,]$.
By then letting $x$ take integral values,
he was able to derive expansions such as the following:
{\allowdisplaybreaks
\begin{multline*}
 \sum\sb {n \geq 0} \frac{1 }{T\sb {4\sp n} (2)}
= [0;1,1,23,1,2,1,18\,815,3,1,23,3,1,23,1,2,1,\\
106597754640383,3,1,23,1,3,23,1,3,18815,1,2,1,23,3,1,23,
 \cdots ],
\end{multline*}
}
where $T_{l}(x)$ denotes the $l$-th Chebyshev polynomial,
and also to derive the continued fraction expansion for
certain sums of series.

At the end of Cohn's paper he listed a number of open questions
and conjectures. One of the problems he mentioned was finding a similar
classification of all those polynomials $f(x) \in \mathbb{Z}[\,x\,]$ for
which the regular continued fraction expansion of the infinite product
{\allowdisplaybreaks
\begin{align}\label{prod1}
\prod_{k=0}^{\infty}\left( 1 +\frac{1}{f_{k}(x)}\right)
\end{align}
}
has  all partial quotients  in $\mathbb{Z}[\,x\,]$.

This turns out to be a technically more difficult problem.  One
reason is that, given any positive integer $k$,  there are classes of
polynomials such as $f(x,k)=2x + x^{2} + x^{k}( (-1)^{k} +(1+
x)g(x))$ for which the regular continued fraction expansion of the
product $\prod_{n=0}^{k}\left ( 1 + 1/f_{n}(x) \right )$ is
specializable for all polynomials $g(x)$. However,    the
specializability/non-specializability of the regular continued fraction expansion of
$\prod_{n=0}^{k+1}\left ( 1 + 1/f_{n}(x) \right )$ does not
immediately follow. This is in contrast to the infinite series
case dealt with by Cohn, where $\sum_{k=0}^{\infty}1/f_{k}(x)$ had
a specializable continued fraction expansion if and only if
$\sum_{k=0}^{2}1/f_{k}(x)$ had a specializable continued fraction
expansion.

In this paper we give several infinite classes of polynomials for which
 $\prod_{n=0}^{\infty}\left ( 1 + 1/f_{n}(x) \right )$ has a specializable
regular continued fraction.

For the polynomials in these classes of degree at least three, we specialize the product at
\eqref{prod1} by letting $x$ take positive integral values,
 producing certain classes of real numbers.  We examine the corresponding
regular continued fractions to  prove the transcendence of these
numbers.

{\allowdisplaybreaks
\section{Some Preliminary Lemmas } \label{prelems}
}

Notation: Unless otherwise stated $f(x)$, $G(x)$, $g(x)$
 will denote  polynomials in
$\mathbb{Z}[\,x ]$, $f_{0}(x) := x$ and, for $n \geq 0, \,\,
f_{n+1}(x) := f(f_{n}(x))$. Sometimes, for clarity and if there is no
danger of ambiguity, $f(x)$ will be
written as $f$ and   $f_{n}(x)$ as $f_{n}$.
 $(f(x))^{m}$ will be written as
$f^{m}$, $(f_{n}(x))^{m}$ as $f_{n}^{m}$ etc.

We will write
\[
\prod_{n}(f(x))=\prod_{n}(f)=
 \prod_{n}  :=\prod_{i=0}^{n}\left(1+\frac{1}{f_{i}}\right)
\]
 and
\[
\prod_{\infty}(f(x))=\prod_{\infty}(f)=\prod_{\infty}:=\prod_{i=0}^{\infty}\left(1+\frac{1}{f_{i}}\right).
\]
 Similarly,
 $S_{n}(f(x))=S_{n}(f)=S_{n}$ will denote the regular continued
fraction expansion (via the Euclidean algorithm) of $\prod_{n}$
and $S_{\infty}(f(x))=S_{\infty}(f)=S_{\infty}$ will denote the regular continued
fraction expansion of $\prod_{\infty}$ (We will use the more concise forms when
there is no danger of ambiguity).

Unless stated otherwise,
the sequence of partial quotients in $S_{n}$ will be denoted by
$\vec{w_{n}}$, so that $S_{n} = [ \vec{w_{n}} ]$.

If a partial quotient in a continued fraction is a polynomial
in $\mathbb{Z}[x]$, it is said to be \emph{specializable}.
A continued fraction  all of  whose partial quotients are specializable
 is also called specializable. We say that a continued
 fraction $[a_{0},a_{1}, \cdots, a_{n}]$ has \emph{even}
(resp. \emph{odd} ) length if $n$ is even (resp. odd).

Since a form of the folding lemma will be used later, we state and
prove this for the sake of completeness.
In what follows let $\vec{w}$ denote the
word $a_{1},\cdots, a_{n}$, $\overset{\leftarrow}{w}$  the
word $a_{n},\cdots, a_{1}$ and $-\overset{\leftarrow}{w}$  the
word $-a_{n},\cdots, -a_{1}$. Let $A_{i}/B_{i}$ denote the $i$-th convergent
to the continued fraction
$[a_{0},a_{1},\cdots]$. Recall that
{\allowdisplaybreaks
\begin{equation}\label{E:cond1}
A_{n}B_{n-1}-A_{n-1}B_{n} =(-1)^{n-1}.
\end{equation}
}

\begin{lemma}\label{l:fl}\cite{VDPS92}
\[[a_{0};\vec{w},Y,-\overset{\leftarrow}{w}]
= \frac{A_{n}}{B_{n}}\left(1+\frac{(-1)^n}{YA_{n}B_{n}} \right) .\]
\end{lemma}

\begin{proof}
It can easily be shown by induction that
{\allowdisplaybreaks
\begin{align*}
[-a_{n};\,-a_{n-1},\cdots, -a_{1}] &= -B_{n}/B_{n-1}\\ \Rightarrow
[a_{0};\vec{w},Y,-\overset{\leftarrow}{w}] &=
[a_{0},\vec{w},Y,-B_{n}/B_{n-1}]
=[a_{0};\vec{w},Y-B_{n-1}/B_{n}]\\ &=
\frac{A_{n}(Y-B_{n-1}/B_{n})+A_{n-1}}
{B_{n}(Y-B_{n-1}/B_{n})+B_{n-1}}\\
&= \frac{A_{n}}{B_{n}}\left(1+\frac{(-1)^n}{YA_{n}B_{n}} \right).
\end{align*}
}
\end{proof}
The last equality uses equation \eqref{E:cond1}.

There are other forms of  symmetry which will appear
later so we give the lemma below. Note that in all of these
 cases $a_{0} = 1$. We  call these symmetries "doubling"
symmetries, following Cohn \cite{C96}.

\begin{lemma}\label{l:fl2}
{\allowdisplaybreaks
\begin{equation}\label{l:fl2a}
[1;\vec{w},Y,-\vec{w}] =
 \frac{A_{n}}{B_{n}}\left(1+\frac{(-1)^n}{A_{n}(B_{n}(Y+1)-A_{n}+B_{n-1})}
 \right).
 \end{equation}
}
{\allowdisplaybreaks
\begin{equation}\label{l:fl2b}
[1;\vec{w},Y,-\overset{\leftarrow}{w},-1]
= \frac{A_{n}}{B_{n}}\left(1+\frac{1}{(-1)^nYA_{n}B_{n}-1} \right).
 \end{equation}
}
{\allowdisplaybreaks
\begin{equation}\label{l:fl2c}
[1;\vec{w},Y,\overset{\leftarrow}{w},1]
= \frac{A_{n}}{B_{n}}\left(1+\frac{1}
{(-1)^{n}B_{n}(Y A_{n}+2A_{n-1})-1} \right).
 \end{equation}
}
{\allowdisplaybreaks
\begin{equation}\label{l:fl2d}
[1;\vec{w},Y,\vec{w}] =
 \frac{A_{n}}{B_{n}}\left(1+\frac{(-1)^n}{A_{n}(B_{n}(Y-1)+A_{n}+B_{n-1})}
 \right).
 \end{equation}
}
\end{lemma}

\begin{proof}
These all follow easily from
 the facts that, for $j = 0,1$,
{\allowdisplaybreaks
\begin{align*}
[(-1)^{j}\vec{w}]&=(-1)^{j}\frac{B_{n}}{A_{n}-B_{n}},\\
[(-1)^{j}\overset{\leftarrow}{w},(-1)^{j}]&=(-1)^{j}\frac{A_{n}}{A_{n-1}}.
\end{align*}
}
We give the proof only for \eqref{l:fl2a} as \eqref{l:fl2b}  \eqref{l:fl2c}
and  \eqref{l:fl2d} follow similarly.
\[
[1;\vec{w},Y,-\vec{w}]
= \frac{A_{n}\left(Y+1 -\frac{A_{n}}{B_{n}}\right)+A_{n-1}}
    {B_{n}\left(Y+1 -\frac{A_{n}}{B_{n}}\right)+B_{n-1}}.
\]
The result follows from \eqref{E:cond1},
after some simple algebraic manipulation.

\end{proof}
 Cohn proved a version of \eqref{l:fl2c} in \cite{C96}.
We also point out that
the doubling symmetry described at  \eqref{l:fl2a} occurs with some classes of polynomials,
but we have not found $S_{\infty}$ to be specializable for polynomials in any of these classes.

For future reference we show how the various forms of
symmetry found in the above lemma  will be used.
Suppose that $\prod_{m}$, when expanded as a continued
fraction, is equal to  $S_{m} = [1;\vec{w}]$, that the
ultimate numerator convergent of $S_{m}$ is $p_{m}$ and
the ultimate denominator convergent is $q_{m}$  and that
$p_{m}'$ and $q_{m}'$ are the penultimate numerator and denominator
convergents respectively,
 that $S_{m}$ is specializable and that $S_{m+1}$ is
related to $S_{m}$ in one of the ways shown in Lemma \ref{l:fl}
or Lemma \ref{l:fl2}.
($Y_{m}$ is used here instead of $Y$ to show the dependence on
$m$). Then
{\allowdisplaybreaks
\begin{align*}
\prod_{m+1} = \prod_{m}\left(1+\frac{1}{f_{m+1}}\right)
=\frac{p_{m}}{q_{m}}\left(1+\frac{1}{f_{m+1}}\right).
\end{align*}
}
On the other hand, from the above lemma,
{\allowdisplaybreaks
\begin{align*}
S_{m+1}
= \frac{p_{m}}{q_{m}}
\left(1+\frac{1}{H(p_{m},q_{m},p_{m}',q_{m}',Y_{m})}\right),
\end{align*}
}
where $H(p_{m},q_{m},p_{m}',q_{m}',Y_{m})$ is a polynomial in its
 variables with integral coefficients that is linear in $Y_{m}$.

If solving the equation $f_{m+1}=H(p_{m},q_{m},p_{m}',q_{m}',Y_{m})$
in terms of $Y_{m}$ leaves $Y_{m}$ in $\mathbb{Z}[ x ]$ for all
$m$ then $S_{m}$ is specializable for all $m$.

For later use we also note that if $x|(f+1)$ then $\prod_{m}$
simplifies to leave $f_{m}$ in the denominator and, say, $r_{m}$ in
the numerator. If ($f_{m}, r_{m})=1$ then, up to sign, the final
numerator convergent of $S_{m}$ is $r_{m}$ and the final denominator
convergent is $f_{m}$. A similar situation also holds if
$(x+1)|f$.

As a result of the following lemma, polynomials of degree $2$ and those of
degree $3$ or more will be considered separately.
\begin{lemma}\label{l3}
If $f(x)$ has degree $\geq 3$ then $S_{n+1}$ contains $S_{n}$ at
the beginning of the expansion.
\end{lemma}

\begin{proof}
Suppose $S_{n}= [1;a_{1},\cdots , a_{m}] = p/q$ where the
$a_{i}$'s, $p$ and $q$ are polynomials in $\mathbb{Q}[x]$.
Let $[1;a_{1},\cdots , a_{i}] =: p_{i}/q_{i}$ and suppose that, via
the Euclidean algorithm:
{\allowdisplaybreaks
\begin{align}\label{E:eu}
p&=q+r_{1}.\\
q&=a_{1}r_{1}+r_{2}. \notag \\
r_{1}&=a_{2}r_{2}+r_{3}. \notag \\
&\phantom{a}\vdots  \notag \\
r_{m-2}&=a_{m-1}r_{m-1}+r_{m}. \notag \\
r_{m-1}&=a_{m}r_{m}.  \notag
\end{align}
}
$S_{n+1} = \displaystyle{p/q\left(1+1/f_{n+1}\right)=
p(f_{n+1}+1)/(qf_{n+1})}$
 and to develop the continued fraction expansion of $S_{n+1}$ one can apply
the Euclidean algorithm to this quotient. From \eqref{E:eu}:
{\allowdisplaybreaks
\begin{align*}
p(f_{n+1}+1)&=q\,f_{n+1}+(r_{1}f_{n+1}+p)\\
qf_{n+1}&=a_{1}(r_{1}f_{n+1}+p)+(r_{2}f_{n+1}-a_{1}p)\\
(r_{1}f_{n+1}+p)&=a_{2}(r_{2}f_{n+1}-a_{1}p)+(r_{3}f_{n+1}+p(1+a_{1}a_{2}))\\
&\vdots
\end{align*}
}
Claim:  Let $r_{-1}'=p(f_{n+1}+1)$, $r_{0}'=qf_{n+1}$ and
for $1 \leq i \leq m$, set
\[
 r_{i}'=r_{i}f_{n+1}+(-1)^{i+1}p\,q_{i-1}.
\]
Then,  for $0 \leq i \leq m-1$,
\[
r_{i-1}'=a_{i}r_{i}'+r_{i+1}'.
\]

Proof of Claim: The claim is clearly true for $i=0,1$ ($a_{0}=1$).
From~\eqref{E:eu}, $r_{i+1}=r_{i-1}-a_{i}r_{i}$ and from the recurrence
relation for the $q_{i}$'s,
$q_{i+1}=a_{i+1}q_{i}+q_{i-1}$.
Suppose the claim is true for
$i=0,1,\cdots,j-1$.
{\allowdisplaybreaks
\begin{align*}
r_{j-1}'-a_{j}r_{j}'&=r_{j-1}f_{n+1}+(-1)^{j}p\,q_{j-2}
-a_{j}(r_{j}f_{n+1}+(-1)^{j+1}p\,q_{j-1}).\\
&=(r_{j-1}-a_{j}r_{j})f_{n+1}+(-1)^{j+2}p(q_{j-2}+a_{j}q_{j-1}).\\
&=r_{j+1}f_{n+1}+(-1)^{j+2}p\,q_{j}.\\
&=r_{j+1}'.
\end{align*}
}
Thus the claim is true.
All that remains to prove the lemma is to show that the degree of
$r_{i+1}'$ is less than the degree of $r_{i}'$ for $0 \leq i \leq
m-1$.

 Let the degree of a polynomial $b$ be denoted by $deg(b)$.
From the Euclidean algorithm it follows that $deg(r_{i+1})<deg(r_{i})$.
Suppose $f$ has degree $r\geq 3$ so that $f_{i}$ has degree $r^{i}$ and   thus,
since $\prod_{i=0}^{n}(1+1/f_{i}(x)) = p/q$, that
\[
deg(p),\, deg(q)  \leq 1 + r + r^2 + \cdots r^{n}
= (r^{n+1}-1)/(r-1).
\]
Thus, for $0 \leq i \leq m$,
\[
deg(p\,q_{i}) \leq deg(p\,q)
\leq 2(r^{n+1}-1)/(r-1) < r^{n+1} = deg(f_{n+1}),
\]
since $r \geq 3$. This implies that, for $0 \leq i \leq m-1$,
{\allowdisplaybreaks
\begin{align*}
deg(r_{i+1}')&=deg(r_{i+1}f_{n+1}+(-1)^{i+2}p\,q_{i})
=deg(r_{i+1}f_{n+1})\\
&<deg(r_{i}f_{n+1})=deg(r_{i}f_{n+1}+(-1)^{i+1}p\,q_{i-1}) =deg(r_{i}').
\end{align*}
}
The result follows.
\end{proof}
Note that if $deg(f) = 2$ (so that $deg(f_{j}) = 2^{j}$)
then  the situation can be quite different.
With the notation of Lemma \ref{l3} and its proof,
the expansion $[1;a_{1}, a_{2}, \dots,a_{k}]$ will  be part of the regular continued fraction
expansion of  $\prod_{i=0}^{n+1}(1+1/f_{i}(x))$
only if $deg(r_{i+1}')<deg(r_{i}')$ for $0 \leq i \leq k$. To ensure this we need
\begin{align}\label{deg2eq}
&deg(p\,q_{i}) < deg(r_{i+1}f_{n+1}),& & 0 \leq i \leq k. &
\end{align}

 With no cancellation
in $\prod_{i=0}^{n}(1+1/f_{i}(x))$, $deg(p) = 2^{n+1}-1$. Thus we can be assured that
 $[1;a_{1}, a_{2}, \dots,a_{k}]$ will  be part of the regular continued fraction
expansion of  $\prod_{i=0}^{n+1}(1+1/f_{i}(x))$ if, for $0 \leq i \leq k$,
\begin{equation*}
deg(p\,q_{i}) \leq 2^{n+1}-1 + deg(q_{i}) <deg(r_{i+1}f_{n+1})=deg(r_{i+1}) + 2^{n+1},
\end{equation*}
or
\begin{equation}\label{qreq}
 deg(q_{i}) \leq deg(r_{i+1}).
 \end{equation}
Since $[1;a_{1},\cdots , a_{i}] =p_{i}/q_{i}$, we have that
$deg(q_{i}) = \sum_{i=1}^{k}deg(a_{i})$.
It is  clear from the recurrence relations at
\eqref{E:eu} that $deg(r_{i}) = deg(a_{i+1})+deg(r_{i+1})$. This implies that
\[
deg(r_{i+1}) = deg(q)-\sum_{j=1}^{i+1}deg(a_{j}).
\]
Thus \eqref{qreq} will hold if
\begin{equation}\label{qraeq}
a_{i+1}+2\sum_{j=1}^{i}deg(a_{j}) \leq deq(q) \leq 2^{n+1}-1.
\end{equation}
Finally, we have that
{\allowdisplaybreaks
\begin{equation}\label{deg2con}
\sum_{j=1}^{k+1}deg(a_{j})  \leq 2^{n} \Longrightarrow
S_{n+1} \text{ starts with } [1;a_{1}, a_{2}, \dots,a_{k}].
\end{equation}
}

We return to the case $deg(f) \geq 3$.
The implication of Lemma \ref{l3}, and the remarks following, is that if $deg(f)\geq 2$,
then it makes sense to talk of the
continued fraction expansion of
$\prod_{i=0}^{\infty}\left(1+ 1/f_{i}\right)$ and, furthermore, that
if $deg(f)\geq 3$,
then  $S_{\infty}$ is a specializable continued fraction if and only if
$S_{n}$ is a  specializable continued fraction for each
integer $n \geq 0$.

Remark: At this stage we are not concerned with whether the polynomials
which are the partial quotients in $S_{\infty}$ have negative leading coefficients or
take non-positive values for certain positive integral $x$. Negatives and zeroes are
easily removed from regular continued fraction expansions (see \cite{VDP94}, for
example).

The following lemma means that we get the proof of the specializability of the
regular continued fraction expansion of $\prod_{k=0}^{\infty}(1+1/f_{k}(x))$ for
some classes of polynomials $f(x)$ for free.
\begin{lemma}\label{lem3}
Suppose $S_{\infty}(f)$ is specializable. Define $g(x)$ by
\begin{equation}
g(x)=-f(-x-1)-1.
\end{equation}
Then $S_{\infty}(g)$ is specializable.
\end{lemma}

\begin{proof}
If $\prod_{k=0}^{\infty}(1+1/f_{k}(x))$ has a specializable continued fraction expansion
$S_{\infty}(f(x)) :=[1;a_{1}(x),a_{2}(x),\dots]$, then
$\prod_{k=0}^{\infty}(1+1/f_{k}(-x-1))$ has the specializable continued fraction expansion
\[
S_{\infty}(f(-x-1))=[1;a_{1}(-x-1),a_{2}(-x-1),\dots].
\]
 Let $g(x)$ be defined as in the statement of the lemma.
For $k\geq 0$,
\[
g_{k}(x) = - f_{k}(-x-1)-1.
\]
This is clearly true for $k=0,1$. Suppose it is true for $k=0,1,\dots , m$.
\begin{align*}
g_{m+1}(x)&=g(g_{m}(x))=g(-f_{m}(-x-1)-1)\\
&=-f(-(-f_{m}(-x-1)-1)-1)-1=-f_{m+1}(-x-1)-1.
\end{align*}
Next,
\begin{align*}
\prod_{\infty}(g(x))&=\prod_{k=0}^{\infty}\left ( \frac{1+g_{k}(x)}{g_{k}(x)} \right)
= \prod_{k=0}^{\infty}\left ( \frac{-f(-x-1)}{-f_{k}(-x-1)-1} \right) \\&=
 \prod_{k=0}^{\infty}\left ( \frac{f(-x-1)}{f_{k}(-x-1)+1} \right).
\end{align*}
From what has been said above, the final product has the regular continued fraction expansion
$[0;1,a_{1}(-x-1),a_{2}(-x-1),\dots]$ and is thus specializable.
\end{proof}

We next demonstrate one of the difficulties in trying to arrive at a complete classification of
all polynomials $f(x)$ for which $S_{\infty}(f)$ is specializable.
\begin{proposition}\label{prop1}
Let $k\geq 2$ be an integer and let $g(x) \in \mathbb{Z}[x]$ be such that
$g(x)$ is not the zero
polynomial if $k=2$. Define
\begin{equation}\label{feqa}
f(x)=2x+x^{2}+x^{k}((-1)^{k}+(x+1)g(x)).
\end{equation}
Then $S_{m}(f)$ is specializable for $m \leq k$.
\end{proposition}
\begin{proof}
We need the following identity. Let $m$ be a non-negative integer. Then
\begin{equation}\label{eq1}
\sum_{j=0}^{m}(-1)^{j-1}\binom{j}{m-j}2^{2j-m}=(-1)^{m+1}(m+1).
\end{equation}
This is most easily proved using the WZ algorithm \cite{PWZ96}. We also
need the fact that, for $t$ a non-negative integer,
\begin{equation}\label{eq2}
(1+k)\sum_{m=0}^{t}(-1)^{m+1}(2k+k^{2})^{m} =  k^{t+1}h_{t}(k)+\sum_{m=0}^{t}(-1)^{m+1}k^{m},
\end{equation}
where $h_{t}(k) \in \mathbb{Z}[k]$. This follows upon expanding the left
side of \eqref{eq2} and noting that, for $1 \leq m \leq t$, that the coefficient of $k^{m}$ is
\[
\sum_{j=0}^{m}(-1)^{j-1}\binom{j}{m-j}2^{2j-m} + \sum_{j=0}^{m-1}(-1)^{j-1}\binom{j}{m-1-j}2^{2j-m+1},
\]
and then using \eqref{eq1}.

Next, we will use the doubling symmetry at \eqref{l:fl2a} to
develop the continued fraction expansion of $\prod_{j=0}^{n}(1+1/f_{j})$,
$1 \leq n \leq k$. Note that
$S_{0}=[1;x]$ and $S_{1}=[1;x,-f/(x(x+1)),-x]$.
For $0\leq n \leq k-1$, let $S_{n}$ be denoted $[1;\vec{w_{n}}]$. Let $A_{n}/B_{n}$ denote
the final approximant and $A_{n}'/B_{n}'$ the penultimate approximant of $S_{n}$.
Note also that $x(x+1)|f$, so that
\begin{equation*}
\prod_{j=0}^{n}\left (1+\frac{1}{f_{j}} \right ) =
\displaystyle{
 \frac{f_{n}+1}{\displaystyle{
x\prod_{j=1}^{n}\frac{f_{j}}{f_{j-1}+1} }}}.
\end{equation*}
Each factor in the denominator of the expression on the right divides $f_{n}$, so
that the numerator and denominator are relatively prime and thus, up to sign,
\begin{align}\label{abeq}
&A_{n}=f_{n}+1,&
&B_{n}=x\prod_{j=1}^{n}\frac{f_{j}}{f_{j-1}+1}=B_{n-1}\frac{f_{n}}{f_{n-1}+1}.&
\end{align}
We will see that the "$+$" sign is the correct choice, so that
above equations are true. By
the correspondence between continued fractions and matrices,
{\allowdisplaybreaks
\begin{equation*}
[1;\vec{w_{n}}]
\longleftrightarrow
\left (
\begin{matrix}
A_{n}& A_{n}' \\
B_{n}& B_{n}'
\end{matrix}
\right ).
\end{equation*}
}
We will see that each $\vec{w_{n}}$ has odd length so that
{\allowdisplaybreaks
\begin{align*}
&[1;\vec{w_{n}},Y_{n},-\vec{w_{n}}]
\longleftrightarrow
\left (
\begin{matrix}
A_{n}& A_{n}' \\
B_{n}& B_{n}'
\end{matrix}
\right )
\left (
\begin{matrix}
Y_{n}& 1\\
1&0
\end{matrix}
\right )
\left (
\begin{matrix}
B_{n}& B_{n}' \\
A_{n}-B_{n}& B_{n}'-A_{n}'
\end{matrix}
\right )\\
&=\left (
\begin{matrix}
A_n^2 - A_nB_n(1 + Y_{n}) - B_nA_n'
& -A_nA_n'  + A_nB_n'(1 + Y_{n}) + A_n'B_n' \\
B_n( A_n - B_n - Y_{n}B_n - B_n' )
& -B_nA_n'  + B_nB_n'(1 + Y_{n}) + B_n'^2
\end{matrix}
\right )\\
&=\left (
\begin{matrix}
A_{n+1}& A_{n+1}' \\
B_{n+1}& B_{n+1}'
\end{matrix}
\right )
\longleftrightarrow
[1;\vec{w_{n+1}}] .
\end{align*}
}
If we set
{\allowdisplaybreaks
\begin{equation}\label{yeq}
Y_{n}=-1
- \frac{
- {\left( 1 + f_{n} \right) }^2 +
f_{1 + n} + B_{n}'(1 + f_{n})}
{B_{n}\,\left( 1 + f_{n} \right) },
\end{equation}
}
and assume that $A_{n} =f_{n}+1$ (we already know this up to sign)
and that $\vec{w_{n}}$ has odd length
 (so that $A_{n}' = (-1+B'_{n}A_{n})/B_{n}=(-1+B'_{n}(f_{n}+1))/B_{n}$, by the
determinant formula), then we get
{\allowdisplaybreaks
\begin{equation}\label{mateq}
\left (
\begin{matrix}
A_{n+1}& A_{n+1}' \\
B_{n+1}& B_{n+1}'
\end{matrix}
\right )
=
\displaystyle{
\left (
\begin{matrix}
\displaystyle{
1 + f_{1 + n}}&
\displaystyle{
\frac{1 + f_{n} - B_{n}' - f_{1 + n}\,B_{n}'}{B_{n}}}\\
\displaystyle{
\frac{{f_{1 + n}}}{1 + {f_n}}B_{n}\,}&
\displaystyle{
1- \frac{ f_{1 + n}\,}{1 + f_{n}}{B_n}' }
\end{matrix}
\right ).
}
\end{equation}
}
We first note that
\[
S_{0}(f) =[1;\vec{w_{0}}] = [1;x]
\longleftrightarrow
\left (
\begin{matrix}
x+1& 1\\
x& 1\end{matrix}
\right ),
\]
so that the $n=0$ case of \eqref{yeq} gives
\[
Y_{0} = -\frac{f}{x(x+1)} \in \mathbb{Z}[x].
\]
Note also that $\vec{w_{0}}$ has odd length, and that
$\vec{w_{n+1}}=\vec{w_{n}},Y_{n},-\vec{w_{n}}$ has odd length
for $n \geq 0$. Further, if $A_{n}=f_{n}+1$
and $B_{n}=x\prod_{j=1}^{n}f_{j}/(f_{j-1}+1)$,
then $A_{n+1}=f_{n+1}+1$ and $B_{n+1}=x\prod_{j=1}^{n+1}f_{j}/(f_{j-1}+1)$.
This means that if the regular continued fraction
expansion of $\prod_{j=0}^{n}(1+1/f_{j})$, namely, $[1; \vec{w_{n}}]$  is specializable,
then  $S_{n+1}(f)=[1; \vec{w_{n}},Y_{n},-\vec{w_{n}} ]$
gives the regular continued fraction expansion
of $\prod_{j=0}^{n+1}(1+1/f_{j})$ and that $S_{n+1}(f)$ is specializable, provided
$Y_{n} \in \mathbb{Z}[x]$, where $Y_{n}$ is defined by \eqref{yeq}. This we now
prove. Suppose we have $Y_{i} \in \mathbb{Z}[x]$ for $0 \leq i \leq n-1 \leq k-2$.
We show $Y_{n} \in \mathbb{Z}[x]$

From the definition of $f(x)$ we have that
\[
f_{n+1}= 2f_{n}+f_{n}^2+f_{n}^{k}((-1)^{k}+(1+f_{n})g(f_{n}).
\]
From \eqref{abeq} we have, for $0 \leq n \leq k-1$, that
$B_{n}\,|\,f_{n}^{n+1}\,|\,f_{n}^{k}$. Thus the result will follow if we can show that
\begin{equation}\label{bneq}
B_{n}\,\bigg |\,B_{n}'+\frac{(-f_{n})^{k}-1}{f_{n}+1}=B_{n}'+\sum_{j=0}^{k-1}(-1)^{j+1}f_{n}^{j}.
\end{equation}
Here and subsequently we mean division in $\mathbb{Z}[x]$.
We now use \eqref{abeq} again and the fact that
$B_{n}'=1-B_{n-1}f_{n}/(f_{n-1}+1)$ (clear  from \eqref{mateq}) to get that \eqref{bneq}
will follow if
{\allowdisplaybreaks
\begin{align*}
&B_{n-1}\frac{f_{n}}{f_{n-1}+1}\bigg |
-B_{n-1}'\frac{f_{n}}{f_{n-1}+1}-f_{n}\sum_{j=0}^{k-2}(-1)^{j+1}f_{n}^{j},
\end{align*}
}
or
{\allowdisplaybreaks
\begin{align}\label{beq2}
&B_{n-1}\bigg |
B_{n-1}'+\sum_{j=0}^{k-2}(-1)^{j+1}f_{n}^{j}(1+f_{n-1}).
\end{align}
}
We make use of \eqref{abeq} again to get that
$B_{n-1}\,|\,f_{n-1}^{n}\,|\,f_{n-1}^{k-1}\,|\,f_{n-1}^{k}$,
so that \eqref{beq2} will hold if
{\allowdisplaybreaks
\begin{align}\label{beq2a}
&B_{n-1}\bigg |
B_{n-1}'+\sum_{j=0}^{k-2}(-1)^{j+1}(2f_{n-1}+f_{n-1}^{2})^{j}(1+f_{n-1}).
\end{align}
}
By \eqref{eq2},
\[
\sum_{j=0}^{k-2}(-1)^{j+1}(2f_{n-1}+f_{n-1}^{2})^{j}(1+f_{n-1})
=\sum_{j=0}^{k-2}(-1)^{j+1}f_{n-1}^{j} + f_{n-1}^{k-1}h(f_{n-1}),
\]
with $h(z) \in \mathbb{Z}[z]$. Since $B_{n-1}'\,|\,f_{n-1}^{k-1}$, we can ignore the second
term on the right above and increase the index on the sum from $k-2$ to $k-1$ for
free, and get that \eqref{beq2a} will hold if
{\allowdisplaybreaks
\begin{align}\label{beq2b}
&B_{n-1}\bigg |
B_{n-1}'+\sum_{j=0}^{k-1}(-1)^{j+1}f_{n-1}^{j}.
\end{align}
}
However, this is just \eqref{bneq} with $n$ replaced by $n-1$ and is thus true by
induction.
Hence the result.
\end{proof}

Remark: One reason we proved Proposition \ref{prop1}
  was to show  that it is not possible
 to eliminate all classes of polynomials
for which $S_{\infty}$ is not specializable by simply
looking at the continued fraction expansion of a finite number of
terms of the infinite product for a general polynomial(Cohn was able to do this in the infinite
series case by looking at just the first three terms).
The polynomial $g(x)$ is arbitrary and it is not
immediately clear whether or not there are certain integers $k$ and certain classes
of polynomials $g(x)$ for which $S_{\infty}(f)$ is specializable. Of course if we restrict $g(x)$ in
Proposition \ref{prop1} to have the form $g(x) =(-1)^{k+1}+x\,h(x)$, then
$S_{k+1}(f)$ will be specializable, since now $f$ will have the formed defined at
\eqref{feqa}, with $k$ replaced by $k+1$. However, it may be that for some $k$ and some other
class of polynomials $g(x)$, that some other form of duplicating symmetry might come into play,
with the result that $S_{\infty}(f)$ becomes specializable.

We remark that Lemma \ref{lem3} gives that if $f(x) = -x^{2}-(1+x)^{k}(1+(-1)^{k+1}\,x\,g(x))$,
with $g(x) \in \mathbb{Z}[x]$, then $S_{k}(f)$ is specializable.

\section{Specializability of $S_{\infty}$ for various
infinite families of  polynomials of degree greater than two }

We can now show that the specializability of $S_{n}$ occurs for
all $n$ for all polynomials in several infinite families. We have
the following theorem.
\begin{theorem}\label{t1}
Let $f(x)$, $G(x)$ and $g(x)$ denote non-zero polynomials in $\mathbb{Z}[x]$
such that the degree of $f(x)$ is at least three.
If $f(x)$ has one of the following forms,
{\allowdisplaybreaks
\begin{align*}
(i)\, f(x)& = x^{2}(x+1)g(x),\\
(ii)\, f(x) &= x(x+1)G(x)-x-1,\\
(iii)\,f(x)  &= x(x+1)G(x)-1= x(x+1)^{2}g(x)-1\\
(iv)\,f(x) &= x(x+1)G(x)-2x-1= x(x+1)((x-1)g(x)+2)-2x-1\\
(v)\,f(x)&= x(x+1)G(x)-2x-2 = x(x+1)((x+2)g(x)-2)-2x-2\\
(vi)\,f(x)&= x(x+1)G(x)-x = x(x+1)(x(x-1)g(x)+1)-x,\\
(vii)\,f(x) &= x(x+1)((x+2)(x+1)g(x)-1)-x-2,
\end{align*}
}
then, for each $n \geq 0$,  $S_{n}$ is a specializable continued fraction.
Hence $S_{\infty}$ is a specializable continued fraction.
\end{theorem}

\begin{proof}
We note that the proof of (iii) follows from the proof of (i) and Lemma \ref{lem3} and
that the proof of (v) likewise follows from the proof of (iv) and Lemma \ref{lem3}.
However, we give an independent proof of (iii) since we also wish to
 demonstrate the type of doubling symmetry exhibited by the
corresponding continued fractions

(i) The product form of the folding lemma will be used. $S_{0} = [1;\,x]$
is clearly specializable.
 Suppose  that
$S_{m}$ is specializable.
For $i \geq 0$
let $S_{i}
=: A_{i}/B_{i}$, where $(A_{i},B_{i})=1$.
From Lemma~\ref{l:fl} it is clear that $S_{m+1}$ is
specializable if
 $A_{m}B_{m}|f_{m+1}$ in $\mathbb{Z}[ x ]$.
After cancellation,
\[
\Pi_{i} = \frac{f_{i}+1}{x\prod_{j=0}^{i-1}f_{j}^{2}g(f_{j})}.
\]
Since $f_{j}|f_{j+1}$ for $j \geq 0$, each term in the denominator of the expression divides
$f_{i}$ and thus the numerator and denominator are relatively prime.
Thus, up to sign, $A_{i} = f_{i}+1$ and
$B_{i}=f_{i-1}^{2}g(f_{i-1})B_{i-1}$. (The first  of these holds  for $i \geq 0$ and the
second for $i \geq 1$).
 It follows easily by induction that $B_{i}|(f_{i})^{2}$.
\[
B_{m}|(f_{m})^{2} \Rightarrow A_{m}B_{m}
=\pm(f_{m}+1)B_{m}|(f_{m}+1)(f_{m})^{2}\Rightarrow
A_{m}B_{m}|f_{m+1}.
\]
Hence the result.

\vspace{15pt}

(ii) This will be proved by inductively constructing the partial quotients
 of $S_{m}$ for each $m$, showing that these partial quotients are in
$\mathbb{Z}[ x ]$ and that the final approximant of $S_{m}$ is equal to
$\prod_{m}$. For this class of polynomial it happens that $S_{m}$
is derived from $S_{m-1}$ by adding a single new partial quotient. We first
need some preliminary results.

For each $i \geq 0$ let $S_{i}=p_{i}/q_{i}$, where $p_{i}$ and $q_{i}$
are the final numerator- and denominator convergents.
It is clear from the definition of $f(x)$ that, for $i \geq 0$,
{\allowdisplaybreaks
\begin{align}\label{E:c4a}
&(f_{i}+1)|f_{i+1},& &f_{i}|(f_{i+1}+1),& &f_{i}|f_{i+2}.&
\end{align}
}
This implies that
{\allowdisplaybreaks
\begin{equation}\label{E:c4ba}
\Pi_{i} = \frac{f_{i}+1}{x\prod_{j=0}^{i-1}(f_{j}G(f_{j})-1)}
=\frac{(x+1)\prod_{j=0}^{i-1}((f_{j}+1)G(f_{j})-1)}{f_{i}}.
\end{equation}
}
This gives  that
$p_{i}|(f_{i}+1)$, and $q_{i}|f_{i}$, for all $i \geq 0$. Next,
\begin{align*}
\frac{p_{i+2}}{q_{i+2}}&=
\frac{p_{i}}{q_{i}}\frac{(f_{i+1}+1)}{f_{i+1}}\frac{(f_{i+2}+1)}{f_{i+2}}
=\frac{p_{i}}{q_{i}}\frac{((f_{i+1}+1)G(f_{i+1})-1)}{(f_{i+1}G(f_{i+1})-1)}
\end{align*}
Claim:
\[
(p_{i},f_{i+1}G(f_{i+1})-1)=
(q_{i},(f_{i+1}+1)G(f_{i+1})-1)=1,
\]
 so that, up to sign,
\begin{align}\label{E:c4b}
p_{i+2}&=((f_{i+1}+1)G(f_{i+1})-1)p_{i} \\
q_{i+2} &= (f_{i+1}G(f_{i+1})-1)q_{i}. \notag
\end{align}
That $(q_{i},(f_{i+1}+1)G(f_{i+1})-1)=1 $ is easily seen to be true since
 $q_{i}|f_{i}$, $f_{i}|f_{i+2}$, so that $q_{i}|f_{i+2}$, but
$((f_{i+1}+1)G(f_{i+1})-1)|(f_{i+2}+1)$. The proof that
$(p_{i},f_{i+1}G(f_{i+1})-1)=1$  is similar.
We are now ready to prove that $S_{n}$ is specializable for $n \geq 0$.

Initially,
$S_{0} = [1,x]$ and $S_{1}=[1;x,-G]$.
It will be shown by induction that
$S_{i}=[\alpha_{0},\alpha_{1},\cdots,\alpha_{i+1}]$, where all the
$\alpha_{j}'s \in \mathbb{Z}[ x ]$ and $(-1)^{i}f_{i}=p_{i-1}q_{i}$.
Both statements are  easily seen to be true for $i=0,1$.
\footnote{We hope our identification $S_{i}=p_{i}/q_{i}$ does not cause
confusion. Usually, if a single continued fraction
$S=[\alpha_{0},\alpha_{1},\cdots,\alpha_{i+1}]$ were being considered,
the value of $S$ would be denoted $p_{i+1}/q_{i+1}$.}

Suppose these statements are true for $i=0,1,\cdots,m-1$.
Let $S_{m-1}=[\alpha_{0},\alpha_{1},\cdots,\alpha_{m}]$.
Set
\begin{equation}\label{aliieq}
\alpha_{m+1} = -\frac{(f_{m-1}+1)}{p_{m-1}}G(f_{m-1})p_{m-2},
\end{equation}
which is in $\mathbb{Z}[ x ]$, since $p_{m-1}|(f_{m-1}+1)$, by the remark
following \eqref{E:c4ba}.
Let $A_{m+1}$ be the final numerator convergent of
$[\alpha_{0},\alpha_{1},\cdots,\alpha_{m},\alpha_{m+1}]$
and let $B_{m+1}$ be its
final denominator convergent.
\begin{align*}
A_{m+1}&=\alpha_{m+1}p_{m-1}+p_{m-2}
= -( (f_{m-1}+1)G(f_{m-1})-1)p_{m-2}\\
B_{m+1}&=\alpha_{m+1}q_{m-1}+q_{m-2}
= -( f_{m-1}G(f_{m-1})-1)q_{m-2}
\end{align*}
The final equality for $B_{m+1}$ uses the facts that
$p_{m-1}q_{m-2}-p_{m-2}q_{m-1}=(-1)^{m-1}$ and
$(-1)^{m-1}f_{m-1}=p_{m-2}q_{m-1}$. Hence, by~\eqref{E:c4b},
$A_{m+1}/B_{m+1}$ $=$  $p_{m}/q_{m}$  $=$ $\prod_{m}$ and
$S_{m}=[\alpha_{0},\alpha_{1},\cdots,\alpha_{m},\alpha_{m+1}]$ . Finally,
\begin{align*}
p_{m-1}q_{m}
 &= p_{m-1}(\alpha_{m+1}q_{m-1}+q_{m-2})\\& =
-(f_{m-1}+1)G(f_{m-1})p_{m-2}q_{m-1} + p_{m-1}q_{m-2}\\
& =
-(f_{m-1}+1)G(f_{m-1})(-1)^{m-1}f_{m-1} + (-1)^{m-1}f_{m-1}
+(-1)^{m-1}\\&
= (-1)^{m}(f_{m-1}+1)(f_{m-1}G(f_{m-1})-1) = (-1)^{m}f_{m}.
\end{align*}
The third equality also uses the facts that
$p_{m-2}q_{m-1} = (-1)^{m-1}f_{m-1}$ and
$p_{m-1}q_{m-2}-p_{m-2}q_{m-1}=(-1)^{m-1}$. Hence $S_{n}$ is
specializable for all $n$.

\vspace{15pt}

(iii) We will use the doubling symmetry found in~\eqref{l:fl2b}.
Suppose $S_{m} = [1,\vec{w_{m}}]$.
It will be shown that $Y_{m}$ can be chosen such that
 \begin{align}\label{eq3}
&\prod_{i=0}^{m+1}\left( 1+\frac{1}{f_{i}} \right ) =S_{m+1}:=
[1,\vec{w_{m}},Y_{m},-\overset{\leftarrow}{w_{m}},-1],& &
Y_{m} \in \mathbb{Z}[ x ].&
\end{align}

Note that $S_{0}=[1;x]$ and that $S_{1}=[1,x,-G,-x,-1]$. $S_{1}$
has even length and if $S_{2}, \dots,  S_{m}$ have been defined
using \eqref{eq3}, then $S_{m}$ has even length. It can be seen
from~\eqref{l:fl2b} that if $S_{m}=p_{m}/q_{m}$ and has even
length, then $f_{m+1}=p_{m}q_{m}Y_{m}-1$ and $Y_{m} \in
\mathbb{Z}[ x ]$ if $p_{m}q_{m}|(f_{m+1}+1)$. This follows easily
by induction.  After cancellation,
\[
\prod_{j=0}^{i}\left( 1+\frac{1}{f_{j}} \right )
=\frac{(x+1)\prod_{j=0}^{i-1}(f_{j}+1)^{2}g(f_{j})}{f_{i}},
\]
so that $p_{i}|(x+1)\prod_{j=0}^{i-1}(f_{j}+1)^{2}g(f_{j})$
and $q_{i}|f_{i}$. Thus it will be sufficient to show that
\begin{align*}
&f_{i}(x+1)\prod_{j=0}^{i-1}(f_{j}+1)^{2}g(f_{j})\,|\,(f_{i+1}+1),& &i\geq 0.&
\end{align*}
This is
clearly true for $i=0$. Suppose it is true for $i=0,1,\dots, m-1$. Then
\begin{align*}
&(x+1)\prod_{j=0}^{m-2}(f_{j}+1)^{2}g(f_{j})\,|\,(f_{m}+1)\\
 \Longrightarrow
f_{m}&(x+1)\prod_{j=0}^{m-1}(f_{j}+1)^{2}g(f_{j})\,| \,
f_{m}(1+f_{m})^{2}g(f_{m})=f_{m+1}+1.
\end{align*}
This completes the proof of (iii).

\vspace{15pt}

(iv) The argument is similar to that used in the proof of (iii).
Here we use the doubling symmetry found in~\eqref{l:fl2c}. Suppose
$S_{m}=[1,\vec{w_{m}}]$ is specializable, has final numerator
convergent $p_{m}$, final denominator convergent  $q_{m}$ and that
$p_{m}'$ and $q_{m}'$ are the penultimate numerator and
denominator convergents respectively.

Note that $S_{1} = [1,x,-G,x,1]$ and
by induction we assume $S_{m}$ has the  symmetric form exhibited in
~\eqref{l:fl2c}, so that $p_{m}'=q_{m}$. Note also that the induction means
that $S_{m}$ has even length, since the duplicating formula always
produces a continued fraction of even length.

It can be seen from~\eqref{l:fl2c} that
\[
S_{m+1}:=[1,\vec{w_{m}},Y_{m},\overset{\leftarrow}{w_{m}},1]
\]
will equal $\Pi_{m+1}$   and be specializable if the equation
$f_{m+1} = q_{m}( p_{m}Y_{m}+2p_{m}')-1$ leads to $Y_{m}  \in
\mathbb{Z}[ x ]$. After cancellation,
\begin{equation}\label{4eq}
\prod_{j=0}^{i}\left( 1+\frac{1}{f_{j}} \right )
=\frac{(x+1)\prod_{j=0}^{i-1}((f_{j}+1)G(f_{j})-2)}{f_{i}}.
\end{equation}
Also,
\[
f-1=(x+1)(x-1)(x\,g(x)+2) \Longrightarrow \prod_{j=0}^{i}(1+f_{j})|(f_{i}^{2}-1),
\]
so that the numerator and denominator in \eqref{4eq} above are
relatively prime. Thus, up to sign $q_{i}=f_{i}$ and $p_{i}\,|\,(f_{i}^{2}-1)$.
\begin{align*}
&f_{m+1}+1=q_{m}( p_{m}Y_{m}+2p_{m}')=q_{m} p_{m}Y_{m}+2q_{m}^{2}
=q_{m} p_{m}Y_{m}+2f_{m}^{2}\\
\Longleftrightarrow & f_{m}(f_{m}^{2}-1)g(f_{m})+2f_{m}^{2}
=q_{m} p_{m}Y_{m}+2f_{m}^{2} \\
\Longleftrightarrow  & Y_{m}= \frac{f_{m}}{q_{m}}\frac{f_{m}^{2}-1}{p_{m}}g(f_{m}).
\end{align*}
Thus $Y_{m} \in \mathbb{Z}[x]$.

Cohn also gave a proof of (iv) in \cite{C96}.

\vspace{15pt}

(v) Here we use the doubling symmetry found at~\eqref{l:fl2d}.
Suppose $S_{m}$ is specializable, has final numerator convergent $p_{m}$,
final denominator convergent  $q_{m}$ and that
$p_{m}'$ and $q_{m}'$ are the penultimate numerator and denominator
convergents respectively. Since $S_{1} = [1,x,-G,x]$ and
$\vec{w_{i}}$  symmetric implies
$ \vec{w_{i}},Y_{i},\vec{w_{i}}$ is symmetric,
we have by induction that $S_{m}$ has odd length and that $\vec{w_{m}}$
is symmetric. This
gives that $q_{m}'=p_{m}-q_{m}$.

It can thus be seen from~\eqref{l:fl2d} that
$[1,\vec{w_{m}},Y_{m},\vec{w_{m}} ]$  will equal $S_{m+1}$ and be
specializable if the equation $f_{m+1} =- p_{m}(
q_{m}(Y_{m}-2)+2p_{m})$ leads to $Y_{m}  \in \mathbb{Z}[ x ]$.
After cancellation, {\allowdisplaybreaks
\begin{equation}\label{5eq}
\prod_{j=0}^{i}\left( 1+\frac{1}{f_{j}} \right )
=\frac{f_{i}+1}{x\prod_{j=0}^{i-1}(f_{j}G(f_{j})-2)}.
\end{equation}
}
Further,
{\allowdisplaybreaks
\begin{equation*}
f+2=(x+2)x((x+1)g(x)-2) \Longrightarrow
\prod_{j=0}^{i}f_{i}|(f_{i+1}+2).
\end{equation*}
}
Since $(f_{i+1}+2) \equiv 2 \mod (f_{i}+1)$ and $f_{i}+1$ is odd for all
$i$ and $x$, the numerator and denominator in \eqref{5eq} above
are relatively prime so that, up to sign,
$p_{i}= f_{i}+1$ and $q_{i}\,|\,(f_{i+1}+2)$.
The result now follows since
{\allowdisplaybreaks
\begin{align*}
f_{m+1} =(f_{m}+1)f_{m}(f_{m}+2)g(f_{m})-2(f_{m}+1)^2&
=- p_{m} q_{m}(Y_{m}-2)-2p_{m}^2 \\
\Longleftrightarrow
(f_{m}+1)f_{m}(f_{m}+2)g(f_{m})&=- p_{m} q_{m}(Y_{m}-2)\\
\Longleftrightarrow Y_{m}-2 &
= \pm \frac{f_{m}(f_{m}+2)}{q_{m}}g(f_{m}).
\end{align*}
}
The last equality follows since $p_{m}=f_{m}+1$, up to sign,
and $f_{m}(f_{m}+2)/q_{m} \in \mathbb{Z}[ x ]$ follows since
$q_{m}\,|\,\prod_{j=0}^{m}f_{j}\,|\,f_{m}(f_{m}+2)$.

\vspace{15pt}

(vi) We first consider the general case where going from $\prod_{i}$ to
$\prod_{i+1}$ adds two terms to the continued fraction expansion,
say $\alpha_{i+1}$ and $\beta_{i+1}$.
Suppose $S_{i}$ is specializable, has final numerator convergent $p_{i}$,
final denominator convergent  $q_{i}$ and that
$p_{i}'$ and $q_{i}'$ are the penultimate numerator and denominator
convergents respectively. Let
$S_{i} = [1,x,\cdots, \alpha_{i},\beta_{i}]$ and
$S_{i+1} = [1,x,\cdots, \alpha_{i},\beta_{i},
\alpha_{i+1},\beta_{i+1}]$. Since $S_{0}=[1,x]$, each $S_{i}$ has odd length.
Using the standard relationship between matrices and continued
fractions, we have that
{\allowdisplaybreaks
\begin{align}
\left(
\begin{matrix}
p_{i+1} & p_{i+1}' \\
q_{i+1} & q_{i+1}'
\end{matrix}
\right)
&=\left(
\begin{matrix}
p_{i} & p_{i}' \\
q_{i} & q_{i}'
\end{matrix}
\right)
\left(
\begin{matrix}
\alpha_{i+1} & 1 \\
1 & 0
\end{matrix}
\right)
\left(
\begin{matrix}
\beta_{i+1} & 1 \\
1 & 0
\end{matrix}
\right)\\
&=\left(
\begin{matrix}
p_{i}(\alpha_{i+1} \beta_{i+1} +1)+p_{i}' \beta_{i+1}
& p_{i} \alpha_{i+1} + p_{i}' \\
q_{i}(\alpha_{i+1} \beta_{i+1} +1)+q_{i}' \beta_{i+1}
& q_{i} \alpha_{i+1} + q_{i}'
\end{matrix}
\right). \notag
\end{align}
}
Using the facts that $S_{i}$ has odd length and that
$S_{i+1} = S_{i}(1+1/f_{i+1})$ it follows that
{\allowdisplaybreaks
\begin{equation}\label{6eq}
f_{i+1} =\frac{p_{i}q_{i}}{-\beta_{i+1}}( \alpha_{i+1} \beta_{i+1}+1)
- p_{i}'q_{i}-1.
\end{equation}
}

We now consider the particular class of polynomials in (vi) above.
Note that
{\allowdisplaybreaks
\begin{align*}
f&=x^2((x^2-1)g(x)+1)& & \Longrightarrow & &
\prod_{j=0}^{i}f_{j}|f_{i+1},\\
f-1& = (x+1)(x-1)(x^{2}g(x)+1) & &
\Longrightarrow & &
\prod_{j=0}^{i}(f_{j}+1)|(f_{i+1}-1),
\end{align*}
}
 so the numerator and denominator of $\Pi_{i}$ are relatively prime and,
up to sign, $p_{i}= \prod_{j=0}^{i}(f_{j}+1)$
and $q_{i}=\prod_{j=0}^{i}f_{j} $. For $i \geq 0$, define
{\allowdisplaybreaks
\begin{align}\label{6abeq}
\alpha_{i+1} &= -\frac{(f_{i}-1)f_{i}g(f_{i})}
{\prod_{j=0}^{i-1}f_{j}(f_{j}+1)}= -\frac{(f_{i+1} -
f_{i}^{2})}{p_{i}q_{i}},\\
\beta_{i+1} &= -\prod_{j=0}^{i}f_{j}(f_{j}+1)=-p_{i}q_{i}. \notag
\end{align}
}
Note that  $\alpha_{i+1}
\in \mathbb{Z}[ x ]$, since $\prod_{j=0}^{i-1}f_{j}(f_{j}+1)\, |\, f_{i}(f_{i}-1)$.
With these values,
{\allowdisplaybreaks
\begin{equation}\label{feq}
f_{i+1}=\alpha_{i+1} \beta_{i+1} -p_{i}'q_{i}
\end{equation}
}
 and
{\allowdisplaybreaks
\begin{align*}
\left(
\begin{matrix}
p_{i+1} & p_{i+1}' \\
q_{i+1} & q_{i+1}'
\end{matrix}
\right)
&=\left(
\begin{matrix}
p_{i}(\alpha_{i+1} \beta_{i+1}+1 -p_{i}'q_{i}) &  p_{i} \alpha_{i+1} + p_{i}'  \\
q_{i}(\alpha_{i+1} \beta_{i+1} -p_{i}'q_{i}) & q_{i} \alpha_{i+1} + q_{i}'
\end{matrix}
\right)\\
&=\left(
\begin{matrix}
p_{i}(f_{i}+1) &  p_{i} \alpha_{i+1} + p_{i}'  \\
q_{i}f_{i}& q_{i} \alpha_{i+1} + q_{i}'
\end{matrix}
\right).
\end{align*}
}
Note that  the
matrix equation above gives that $p_{i}=\prod_{j=0}^{m}(f_{j}+1)$,
$q_{i}=\prod_{j=0}^{m}f_{j}$, for all $i$
(looking at $\prod_{i}$ gave this only up to sign).

Upon expanding $\Pi_{1}$ we have that $S_{1}=[1,x,-G,-x(x+1)]$,
satisfying \eqref{6abeq}. We now show that if the $\alpha_{i}$'s
and $\beta_{i}$'s are defined inductively by \eqref{6abeq}, then
for $i\geq 0$, $S_{i}$ is specializable and $S_{i}=\Pi_{i}$.
Suppose these hold for $i=0,1,\dots, m$.

 On substituting for $\alpha_{m+1}$ and $\beta_{m+1}$
 in \eqref{feq} it follows that
{\allowdisplaybreaks
\begin{equation*}
f_{m+1} = f_{m+1} - f_{m}^{2} - p_{m}'q_{m} \Longrightarrow
p_{m}' =-\frac{f_{m}^{2}}{q_{m}} = -\frac{f_{m}}{q_{m-1}},
\end{equation*}
}
so that all that is necessary to continue the induction is to show
 $p_{m+1}' =-f_{m+1}/q_{m}$.
This follows easily from the
relation (clear from the matrix equation above) that
$p_{m+1}'= p_{m} \alpha_{m+1} + p_{m}'$ and \eqref{6abeq}.
These relations also give that
{\allowdisplaybreaks
\begin{equation*}
q_{m}' = \frac{-f_{m}+1}{p_{m-1}}.
\end{equation*}
}

\vspace{10pt}

(vii) This follows from (vi) and Lemma \ref{lem3}.
\end{proof}

\section{The Degree Two Case}
It is a much easier task to give a complete classification of all
polynomials $f(x)$ of degree two for which $S_{\infty}(f(x))$  is
specializable. Essentially, the method is to start with a general
polynomial $ax^{2}+bx+c$ and to choose an integer $n$ large enough
so that some part of the continued fraction expansion of
$\prod_{k=0}^{n}(1+1/f_{k})$, say $[1,a_{1}(x), \dots , a_{t}(x)]$
forms part of the continued fraction expansion of
$\prod_{k=0}^{\infty}(1+1/f_{k})$ (This follows by virtue of the
remarks following Lemma \ref{l3} for the degree 2 case). The
coefficients  in the $a_{i}(x)$ will be rational functions in $a$,
$b$ and $c$ and the requirement that the $a_{i}(x) \in
\mathbb{Z}[x]$ will impose conditions on $a$, $b$ and $c$, leading
to the following theorem.
{\allowdisplaybreaks
\begin{theorem}\label{t2}
Let $f(x) \in \mathbb{Z}[x]$ be a polynomial of degree two such that $S_{\infty}(f)$
is specializable. Then
{\allowdisplaybreaks
\begin{equation*}
f(x) \in \{-x^2-2x-2, \, -x^2-2x-1,\,x^2, \,x^2-1\}.
\end{equation*}
}
\end{theorem}
}

\section{Specialization and Transcendence}

In what follows, we assume $f(x) \in \mathbb{Z}[x]$ and $M \in \mathbb{Z}$ are such that
$f_{j}(M) \not = 0, -1$, for $j\geq 0$.

For any of the polynomials $f$ in Theorems \ref{t1} and \ref{t2},  $S_{\infty}(f)$
will typically have some partial quotients which are polynomials in $x$ with
leading negative coefficients. It may also happen that if $S_{\infty}(f)$ is specialized
by letting $x$ assume integral values, that negative or zero partial quotients may appear
in the resulting continued fraction.
These are easily removed, as the following equalities show
(see also  \cite{VDP94}).
\begin{align*}
[\dots , a,b,0,c,d,\dots]&=[\dots , a,b+c,d,\dots],\\
[\dots , a,-b,c,d,e,\dots]&=[\dots , a-1,1,b-1,-c,-d,-e,\dots]
\end{align*}
Thus, if $M$ is an integer, repeated application of the identities above
will transform $S_{\infty}(f(M))$ to produce the regular continued fraction expansion of
the corresponding real numbers.

One also easily checks that, in the first equation,  if $P'$ is the numerator
convergent of $[\dots , a,b,0]$ and $P$ is the is the numerator
convergent of $[\dots , a,b,0,c]$, then, on the right side, $P'$ is the numerator
convergent of $[\dots , a]$ and $P$ is the is the numerator
convergent of $[\dots , a,b+c]$. A similar relationship hold for the denominator convergents.
It is clear also that the numerator and denominator convergents agree in both continued fractions
from $d$ onwards.

Similarly, for the continued fraction on the left in
the second equation, if $P'$ is the numerator
convergent of $[\dots , a]$ and $P$ is the is the numerator
convergent of $[\dots , a,-b]$, then, on the right side, $P'$ is the numerator
convergent of $[\dots , a-1,1]$ and $-P$ is the is the numerator
convergent of $[\dots , a-1,1,b-1]$. As above, similar relationship
hold for the denominator convergents. In this second case the numerator and denominator
convergents agree up to sign from $c$ onwards, but only every second one has the same sign.

Suppose, before removing zeroes and negatives, that
{\allowdisplaybreaks
\begin{equation*}
S_{\infty}(f(M))=[1;a_{1},a_{2}, \dots, a_{n},a_{n+1}, a_{n+2},\dots].
 \end{equation*}
}
 Suppose that all the zeroes and negatives are removed from the expansion as far as $a_{n}$,
so that
{\allowdisplaybreaks
\begin{equation*}
S_{\infty}(f(M))=[b_{0};b_{1},b_{2}, \dots, b_{t},\pm a_{n+1}, \pm a_{n+2},\dots],
 \end{equation*}
}
with $b_{0}$ a non-negative integer and $b_{i}$ a positive integer, for $1 \leq i \leq t$.
This transformation will leave
{\allowdisplaybreaks
\begin{equation*}
[1;a_{1},a_{2}, \dots, a_{n}]=[b_{0};b_{1},b_{2}, \dots, b_{t}] := \frac{A_{n}}{B_{n}}.
\end{equation*}
}
How the regular continued fraction expansion of $S_{\infty}(f(M))$ continues will depend on the
signs of $\pm a_{n+1}$ and $\pm a_{n+2}$. We assume $|a_{n+1}|>2$ and $|a_{n+2}|>1$
(This is justified for all the cases we
consider, since $|a_{n+1}|$ will a large positive number and
$|a_{n+2}|$ will be at least $M$). If any negatives in $\pm a_{n+1}$ and $\pm a_{n+2}$
are removed, then we get that
{\allowdisplaybreaks
\begin{multline*}
S_{\infty}(f(M))=\\
\begin{cases}
[b_{0};b_{1},b_{2}, \dots, b_{t}, |a_{n+1}|, \dots], & \pm a_{n+1}>0, \pm a_{n+2}>0,\\
[b_{0};b_{1},b_{2}, \dots, b_{t}, |a_{n+1}|-1, \dots], & \pm a_{n+1}>0, \pm a_{n+2}<0,\\
[b_{0};b_{1},b_{2}, \dots, b_{t}-1, 1,|a_{n+1}|-1, \dots], & \pm a_{n+1}<0, \pm a_{n+2}<0,\\
[b_{0};b_{1},b_{2}, \dots, b_{t}-1,1, |a_{n+1}|-2, \dots], & \pm a_{n+1}<0, \pm a_{n+2}>0.
\end{cases}
\end{multline*}
}
If $b_{t}=1$, then the zero  $b_{t}-1$ in the latter two cases is removed by replacing
the string $b_{t-1},0,1$ by $b_{t-1}+1$. The point of all this is that, using facts from the
theory of regular continued fractions, we can now say that, in all cases,
{\allowdisplaybreaks
\begin{equation}\label{transeq}
\left | \prod_{j=0}^{\infty} \left ( 1+\frac{1}{f_{j}(M)}\right ) - \frac{A_{n}}{B_{n}} \right | =
\left | S_{\infty}(f(M)) - \frac{A_{n}}{B_{n}} \right | < \frac{1}{B_{n}^{2}(|a_{n+1}|-2)}.
\end{equation}
} We are now ready to prove some results on transcendence for the
polynomials $f(x)$ in Theorem \ref{t1}. For simplicity we consider
only the case where each $f(x)$ has positive leading coefficient
and, when specializing the infinite product
$\prod_{j=0}^{\infty}(1+1/f(x))$ and $S_{\infty}(f(x))$ by letting
$x$ be an integer $M$, we assume $M>0$ is large enough so that
$f(M)>0$ and $f(x)$ is strictly increasing for $x \geq M$. We will
use Roth's Theorem.
\begin{theorem}$($Roth \cite{R55A}$)$
Let $\alpha$ be an algebraic number and let $\epsilon>0$. Then the
inequality
{\allowdisplaybreaks
\begin{equation*}
\left | \alpha - \frac{p}{q} \right | <\frac{1}{q^{2+\epsilon}}
\end{equation*}
}
has only finitely many solutions $p \in \mathbb{Z}$, $q \in \mathbb{N}$.
\end{theorem}
 We have the following theorem.
\begin{theorem}\label{t3}
Let $f(x)\in \mathbb{Z}[x]$ be a polynomial with leading positive coefficient
such that $f(x)$ satisfies one of the congruences in Theorem \ref{t1}. Let
$M >0$ be an integer large enough so that $f(M)>0$ and $f(x)$ is strictly increasing
for $x \geq M$. Then the number $\prod_{j=0}^{\infty}(1+1/f(M))$ is transcendental.
\end{theorem}

\begin{proof}
(i) Let $f(x) = x^{2}(x+1)g(x)$ and let
$S_{n}(f(M))=[1;\vec{w_{n}}] = p_{n}/q_{n}$. From Theorem \ref{t1} (i),
\[
S_{n}(f)=[1;\vec{w_{n}},Y_{n},-\overset{\leftarrow}{w}],
\]
with, up to sign, $p_{n}=f_{n}(M)+1$ and $Y_{n} = f_{n+1}/(p_{n}q_{n})$.
Since $p_{n}=f_{n}(M)+1$, $q_{n}$ also has order $f_{n}(M)$. Since $f(x)$ has
degree at least $3$, $Y_{n}$ also has order at least $f_{n}(M)$. Hence, after
any negatives are removed the inequality at \eqref{transeq}, with $a_{n+1}$ being
replaced by $Y_{n}$, will give that the inequality
{\allowdisplaybreaks
\begin{equation*}
\left | S_{\infty}(f(M)) - \frac{p_{n}}{q_{n}} \right | < \frac{1}{q_{n}^{2+\epsilon}}
\end{equation*}
}
will have infinitely many solutions $(p_{n},\,q_{n})$ with $\epsilon = 1/2$, say.
This gives that $S_{\infty}(f(M))$ is transcendental.

\vspace{10pt}

(ii) Let $f(x) = x(x+1)G(x)-x-1$, with $G(x) \in \mathbb{Z}[x]$. We use the notation used
in the proof of (ii), Theorem \ref{t1}. From \eqref{aliieq} it can be seen that
$\alpha_{m+1}$ has order at least $f_{m-1}$, since $G(x)$ has degree at least one.
Once again, after removing negatives, the inequality at \eqref{transeq}, with $a_{n+1}$ being
replaced by $\alpha_{m+1}$, will give that the inequality
{\allowdisplaybreaks
\begin{equation*}
\left | S_{\infty}(f(M)) - \frac{p_{m-1}}{q_{m-1}} \right | < \frac{1}{q_{m-1}^{2+\epsilon}}
\end{equation*}
}
will have infinitely many solutions $(p_{m-1},\,q_{m-1})$ with $\epsilon = 1/2$, say.
Thus  $S_{\infty}(f(M))$ is transcendental.

\vspace{10pt}

(iii), (iv) and (v) The proofs in these case are virtually identical to the proof
for the class of polynomials in (i), and so details are omitted. It suffices to point out,
with the notation of Theorem \ref{t1},  that in each case
 $q_{m}$ has order at most $f_{m}(M)$ and that $Y_{m}$ has order at least $f_{m}(M)$.

\vspace{10pt}

(vi) From \eqref{6eq} and \eqref{6abeq} we have that
{\allowdisplaybreaks
\begin{align*}
|\alpha_{i+1}|
&= \left | g(f_{i}(M)) \prod_{j=0}^{i-1} \{f_{j}(M)^{2}g(f_{j}(M))+1  \}
 \{ (f_{j}(M)^{2}-1)g(f_{j}(M))+1  \} \right |\\
&\geq\left | g(f_{i}(M)) \prod_{j=0}^{i-1}f_{j}(M)^{2}((f_{j}(M)^{2}-1)g(f_{j}(M))+1) \right |\\
& \geq \left | \frac{1}{M} \prod_{j=0}^{i}f_{j}(M)\right | =\frac{ |q_{i}|}{M}.
\end{align*}
}
Thus $\alpha_{i+1}$ has the same order as $q_{i}$ and a similar argument to that
used above applied to \eqref{transeq} with $(A_{n},B_{n})=(p_{n},q_{n})$ and
$a_{n+1}=\alpha_{n+1}$ gives that $\prod_{j=0}^{\infty}(1+1/f_{j}(M))$ is transcendental.

\vspace{10pt}

(vii) We omit the proof for polynomials in this class by an appeal to Lemma \ref{lem3} and
(vi) above. Strictly speaking, to apply Lemma \ref{lem3} we should also have proved
(vi) for polynomials $f(x)$ with negative leading coefficient. An alternative would have been
to give a direct proof of case (vii) in Theorem \ref{t1} and then prove Theorem \ref{t3}
for polynomials in this class in the same way we proved (vi) above,
and this would have not been difficult.

\end{proof}

Remark: In light of Theorem \ref{t3}, an obvious question is the following: If 
$f(x) \in \mathbb{Z}[x]$ is a polynomial of degree at least three 
and $M$ is an integer such that $f_{j}(M) \not =0,-1$ for any $j$ and 
$f_{j}(M) \not = f_{k}(M)$ for $j \not = k$,  is the infinite product 
\[
\prod_{j=0}^{\infty} \left (1 + \frac{1}{f_{j}(M)}\right )
\]
transcendental? If this is false, find a counter-example.

With this question in mind, we note the following "near exception":  if $f(x) = 4 x^{3}+6x^{2}-3/2$
and $M$ is any integer different from $-1$, then 
\[
\prod_{j=0}^{\infty} \left (1 + \frac{1}{f_{j}(M)}\right ) = \sqrt{\frac{2M+3}{2M-1} }.
\]

We now look at some particular examples of specialization. As Cohn showed in \cite{C96}, if
$l \equiv 2 \mod{4}$,  and $T_{k}(x)$ denotes the $k$-th Chebyshev polynomial then
\[
\prod_{j=0}^{\infty} \left ( 1 +\frac{1}{T_{l^{j}}(x)} \right )
\]
has a specializable continued fraction expansion with predictable partial quotients.
This follows from Theorem \ref{t1} (iv), using the facts that $T_{1}(x)=x$,
that if $l \equiv 2 \mod{4}$ then
$T_{l}(x) \equiv 2x^{2}-1 \mod{x(x^{2}-1)}$ and that $T_{a}(T_{b}(x))=T_{ab}(x)$, for
all positive integers $a$ and $b$. For example, setting $l=6$ and $x=3$, we get after
removing negatives, that
\begin{multline*}
\prod_{j=0}^{\infty} \left ( 1 +\frac{1}{T_{6^{j}}(3)} \right ) =\\
 [1, 2, 1, 1632, 1, 2, 1, 3542435884041835200, 1, 2, 1, 1632, 1, 2, 1,\\
 2602
9539217771234538544216588488566196402655804477165253\\9336341
222077618284068468732
496046837200411447595913600,\\
 1, 2, 1, 1632, 1, 2, 1, 3542435884041835200, 1, 2,
1, 1632, 1, 2, 1, \dots ].
\end{multline*}

In part (vi) of Theorem\ref{t1}, setting $g(x) = (x^{2k-2}-1)/(x^{2}-1)$ gives $f(x)=x^{2k}$,
for $k \geq 2$, so that
\[
\prod_{j=0}^{\infty} \left ( 1 +\frac{1}{x^{(2k)^{j}} }\right )
\]
has a specializable continued fraction expansion with predictable partial quotients.
This result can also be found in \cite{MVDPS91}, where the formulae
for the partial quotients that we have are also given.
For example, if $k=2$ and $x\geq 2$ is a
positive integer, then
\begin{align*}
\prod_{j=0}^{\infty} &\left ( 1 +\frac{1}{x^{4^{j}} }\right )=
\big[1;x-1,1,x(x-1),x(x+1), \\
&x^{3}(x-1)(x^{2}+1), x^{5}(x+1)(x^{4}+1),\\
& x^{11}(x-1)(x^{2}+1)(x^{8}+1), x^{21}(x+1)(x^{4}+1)(x^{16
}+1),
\dots,\\
&x^{(2\times 4^{i}+1)/3}(x-1)\prod_{j=0}^{i-1}(x^{2\times 4^{j}}+1),
x^{(4^{i+1}-1)/3}(x+1)\prod_{j=0}^{i}(x^{ 4^{j}}+1), \dots \big ].
\end{align*}

\section{Concluding Remarks}

Ideally, one would like to  have a complete list of all classes of
polynomials $f(x)$  for which $\prod_{n=0}^{\infty}(1+1/f_{n})$
has a specializable continued fraction expansion. We hesitate to
conjecture that our Theorems \ref{t1} and \ref{t2} give such a
complete list, since there may be other classes of polynomials for
which $S_{\infty}$ displays more complicated forms of duplicating
symmetry. One reason for suspecting this is that Cohn \cite{C96}
found some quite complicated duplicating behavior for several
classes of polynomials. One example he gave was the class of
polynomials of the form
\[
f(x)=x^{3}-x^{2}-x+1+x^{2}(x-1)^{2}g(x),
\]
with $g(x) \in \mathbb{Z}[x]$. If $S_{n}=\sum_{j=0}^{n}1/f_{j} = [0;\vec{s_{n}}]$,
then, for $n \geq 3$,
{\allowdisplaybreaks
\begin{equation}\label{deq}
S_{n}=[0;\vec{s_{n-1}},X_{n},-\vec{s_{n-2}},0,\vec{s_{n-4}},Y_{n-2},0,
Z_{n},-\overset{\leftarrow}{s_{n-4}},Y_{n},\overset{\leftarrow}{s_{n-2}}],
\end{equation}
}
where the $X_{i}$, $Y_{i}$ and $Z_{i}$ are polynomials in
$\mathbb{Z}[x]$. It is not unreasonable to suspect similar such
complicated behavior in the infinite product case. It may also happen that
there are classes of polynomials for which a result similar to that in Proposition
\ref{prop1} holds, with the simple duplicating symmetry of Proposition \ref{prop1}
being replaced by some more complicated form of duplication such as that in
Equation \ref{deq} above.

We have already mentioned one difficulty in Proposition \ref{prop1}, where it
was shown that it was not possible to eliminate all classes of polynomials
for which $S_{\infty}$ is not specializable simply by examining $S_{n}$, for
some fixed finite $n$.

Even apart from  this theoretical difficulty, the infinite product
case is computationally more difficult than the infinite series
case. We found that to eliminate some classes of polynomials it
was necessary to go out to $S_{8}$. Our method in this case
involved finding the regular continued fraction expansion of a
rational function of degree more than $16,000,000$ in numerator
and denominator. The rational function had a free parameter and
the class of polynomials was eliminated by showing no value of the
parameter made $S_{8}$ specializable.
 For other classes of polynomial, going out to $S_{8}$ was insufficient.

In conclusion, both theoretical and computational difficulties prevented us from
arriving at a complete classification of all polynomials $f(x)$ for which
$S_{\infty}(f)$ is specializable. However, we hope the results in this paper
will stimulate further work on this problem.

{\allowdisplaybreaks

}
\end{document}